# A bi-objective model to locate several bio-refineries and optimize their supplies

**Corresponding author:**

Nasim ZANDI ATASHBAR
ICD-LOSI, University of Technology of Troyes (UTT)
Troyes, France
CS 42060, 10004
nasim.zandi_atashbar@utt.fr

**Co-Authors:**

Nacima LABADIE
Same address as corresponding author
Nacima.Labadie@utt.fr

Christian PRINS
Same address as corresponding author
Christian.Prins@utt.fr

**Introductive summary**

Many research projects are devoted to biomass production, selection of new crops and bio-refinery processes but in comparison the biomass supply chain to satisfy the planned demands of refineries are relatively neglected, although it could constitute the Achille's heel of the system. Indeed, an important fraction of biomass cost at refinery gates resides in logistic costs. This paper briefly presents a mathematical model to locate several refineries over a large territory and optimize their supplies over one year divided into weeks. Two objectives are minimized: the total cost of the system, GHG emissions, and energy consumptions.

**Purpose**

Growing consciousness about destructive effects of climate change caused by greenhouse gas emissions, along with a huge rise in global demand for energy, have forced many researchers to look for better alternatives to fossil fuels. Biofuel derived from biomass, as a renewable and clean energy source, is one of the few potential replacements of fossil fuels, and can play a crucial role in the transition from traditional sources of energy. Biomass flow from field to fuel is called biomass supply chain and includes various activities such as cultivation, harvesting, handling, storage, transportation, and biofuel conversion. Although biomass itself is cheap relative to other sources of energy, due to high logistics expenses its cost at refinery gates can be decisive. The task of minimizing logistic cost is really challenging, because contrary to industrial logistics, the raw materials (oilseed and lignocellulosic crops) are produced slowly, seasonally, and with a limited yield, over vast territories. In particular, a refinery must use successive crops during the year, e.g., miscanthus in spring, rape in July, cereal straws in August, camelina in

October and short rotation trees like willows in winter. Besides, the economic aspect of biomass supply chain is not the only important aspect. In fact, recently, there has been increasing concerns about environmental aspects as well. Therefore, to make the production of biofuel affordable, it is critical to improve the efficiency of its supply chain by considering economic and environmental aspects simultaneously. Hence, more and more researchers have been involved in modelling and optimizing biomass supply chains [1, 2]. This paper concerns a French national research project on biomass logistics. Its goal is to develop a model to locate several biorefineries on a large territory (equivalent to two Frech administrative regions), and to plan their supplies in biomass over a multiperiod horizon of one year, divided in one-week slots. The criteria to be minimized are the total cost of the system and GHG emissions.

**Approach and scientific innovation**

This paper proposes a multi-objective, multi-period, Mixed Integer Linear Programming (MILP) model to optimize simultaneously the economic and environmental performance of multi-biomass supply chain for several bio-refineries at the tactical level. The first objective is to minimize the total cost of the supply chain, including biomass production, storage, handling, bio-refineries setup and transportation. The second one is to minimize the greenhouse gas (GHG) emissions including biomass productions, farm storages, centralized storages, biofuel production and transportation. The amount of biomass produced, shipped and stored during each period as well as the number, size and locations of bio-refineries are determined. To the best of our knowledge, this is the only research that deals with a multi-objective, multi-period biomass supply chain, at tactical levels, considering different biomass types, centralized storages and by locating new bio-refineries or using the existing refineries.

The MILP is designed to handle comprehensive multi-period and multi-biomass supply chain with several node types. Biomass can be harvested in elementary production zones (small administrative areas called "cantons" in French), and then either stored in farm storages or transferred directly to centralized storages. Biomass can also be shipped from farm storages to centralized storages. Finally, it is transported to the refineries. The supply chain can be described by a graph with a node-set composed of biomass production zones, farm storages, centralized storages and bio-refineries input stocks, and an arc-set. Each arc denotes a pre-computed shortest path between any two nodes in the road network, with specified length and a required vehicle.

The model relies on some assumptions: (1) a refinery will not shut down once it is operational; (2) each refinery is already placed or must be located, and there is at most one per zone; (3) each refinery defines its needs in dry tonnes per product and per period, 4) biomass are transported only by road; (4) the mathematical model is "data driven": all data even the network structure are stored in external files; (5) it is possible to add new products or new facilities; (6) the supply chain ranges from finished products (ready to ship in the farms) to bio-refineries storages; (7) the planning horizons are divided in discrete time slots (period), currently 52 periods of 7 days.

The input data, stored in a data base, include: (1) cost functions associated with production, farm storage, centralized storage, handling, bio-refineries setup and transportation; (2) the geographic distance between each node in the biomass supply chain processed by MapPoint; (3) GHG emissions associated with biomass production zones, farm storages, centralized storages, bio-refineries and transportation; (4) the annual yield of each type of biomass and annual biofuel demand; (5) initial inventory for each node (biomass production zones, farm storages, centralized storages and bio-refineries); (6) Loss coefficient per period for each node; (7) capacities for biomass production zones, farm storages, centralized storages and bio-refineries; (8) harvesting window for each type of biomass.

To minimize the total cost and GHG emissions, the following optimal decisions are determined by the model: selection of biomass production zones; inventory level of harvested biomass in each biomass production zones during time period; material flow from production zones to farm storages; material flow from production zones to centralized storages; Inventory levels of farm storages during time period; material flow from farm storages to centralized storages; Inventory levels of centralized storages during time period; material flow from centralized storages to bio-refineries; Inventory levels of bio-refineries during time period; site selection for location of bio-refineries; capacity level (size) for selected bio-refineries.

As already mentioned, the model is designed to minimize two objectives (optimization criteria). The first objective is to minimize the total cost of biomass at the refinery gates. The total cost includes biomass production costs in farm field, handling costs in the farms, transportation costs from fields to farm storages, farm storage costs, handling costs at farm storages, transportation costs to centralized storages, centralized storage and handling costs, transportation costs from centralized storages to bio-refineries, annualized setup and operating costs of bio-refineries.

The second objective function is to minimize the total greenhouse gas (GHG) emissions, measured as kilograms of equivalent $CO_2$ per tonne. The total GHG emissions include emissions related to biomass productions, farm storages, centralized storages, handling, and transportation. The transportation emissions consist of the different transport steps in the logistic network.

In our model, the constraints are formulated in a generic way to be applied to several types of nodes if possible, although some constraints are specific to the bio-refineries. The storage capacity constraints apply to all production zones, farm storages, centralized storages and bio-refineries. They ensure that each site respects its available capacity. Shared capacity constraints for a group of nodes sharing the same storage capacity for different products, are considered. Also, final inventory in the last period must be respected. Inventory balance constraints control the balances in each node. Two forms for inventory balance are considered and both apply the loss factor to the stock from previous period. The first form is a particular case for the first period, with the initial inventory. Also, the maximum throughput constraints are used to limit the total flow leaving a node. The demand satisfaction constraints are considered and they look like the inventory balance equations, but the output flow is replaced by a demand.

The specific constraints on bio-refinery creations are the following. The user may define several types of refineries (in terms of cost, production capacity, and calendar of demands), a number of refineries for each type, and allowed/forbidden zones to create these refineries. Some refineries can be already located. At most one refinery can be built in each allowed zone and the number of refineries created for each type must be equal to the number specified.

Instead of adding the two objectives, which leads to bad results in case of conflicting objectives, we use Pareto optimization which offers to the decision maker a set of different compromise solutions. The bi-objective problem is solved using the ε-constraint method. The main idea is to transform one objective into a constraint to come back to a classical, single objective problem. Here, the first objective $f(x)$ (the total cost (€), where $x$ is the vector of variables) is minimized, while the minimization of the second objective $g(x)$ (GHG emissions ($CO_2$-eq tons)) is replaced by a constraint $g(x) \leq \varepsilon$. The resulting model is then solved for different values of $\varepsilon$.

**Results**

The model is already tested on realistic data, covering a circle with a 50-km radius centered on the city of Compiègne (60 km north of Paris), with 29 zones (administrative districts containing several communes each, and 1768 farms in total). The farms are located in the departments of Oise (60), Aisne (02) and Somme (80). The 3 rape products (bulk seeds, straw bales and chaff bales), and a 1-year horizon divided into 52 weeks are considered. Rape production in each zone was estimated using results of the 2010 Agricultural Census. Storage capacities and costs for silos (for seeds) and platforms (for bales) were obtained by sending a questionnaire to centralized storage operators. The costs of handling equipment and transport vehicles were found in professional databases. The total demand of biomass is assumed to be 160160 tonne. One refinery is assumed to be already located in Compiègne, while a second may be created in any district with no common border with Compiègne.

The resulting instance was solved using Xpress-IVE 7.8 from FICO, on a 2.70 GHz Intel Core i7 portable PC with 32 GB of RAM and Windows 7 Professional. The model has 108883 variables and 16368 constraints. The pre-solver reduces it to 106582 variables and 11054 constraints.

In this case, the obtained optimum value of objective function (1) is 65 696 648 € in around 49s. The cost of biomass represent 50.13%, capital and operating costs of refineries 30.4%, transport 3.65%, handling 2.12%, and storage 1.7%. The obtained optimum value of second objective function is 68 616 735 $CO_2$-eq tons. By decreasing the value of $\varepsilon$, we observe of course an increase of the costs and a decrease in the amount of GHG emissions, showing the conflicting nature of the two objectives. As expected, the results reveal a trade-off between environmental and economic performance and a curve of compromises can be plotted to help the decision maker.

**Conclusion**

A multi-objective, multi-period, Mixed Integer Linear Programming (MILP) model has been developed to optimize simultaneously the economic and environmental performances of a multi-biomass supply chain for several bio-refineries. The two optimization objectives (total cost and GHG emissions) are minimized using the $\varepsilon$-constraint method. The resolution of the model using a state-of-the art commercial MILP solver shows that large supply chains can be modelled and optimized in acceptable running time on a standard PC.

The work in progress consists in testing the model on larger cases, to see the increase in running time, and to handle multimodal transportation, e.g., using rivers and railways. Concerning the first item, decomposition or relaxation techniques will be applied if the running time becomes excessive. Concerning multimodal transportation, the current modelling of the logistic network is already very generic. It is foreseen to let the user describe several networks (road, rail, river) as graphs, and to interconnect these graphs using transhipment nodes.

**References**
[1] Ba, B. H., Prins, C., & Prodhon, C. (2016). Models for optimization and performance evaluation of biomass supply chains: An Operations Research perspective. Renewable Energy, 87, 977-989.
[2] De Meyer, A., Cattrysse, D., Rasinmäki, J., & Van Orshoven, J. (2014). "Methods to optimise the design and management of biomass-for-bioenergy supply chains: a review." Renewable and sustainable energy reviews, 31, 657-67.